\documentclass[a4paper,10pt]{article}

\def \Z {{\mathbf {Z}}}
\def \R {{\mathbf {R}}}
\def \N {{\mathbf {N}}}
\def \F {{\cal F}}

\def\tr{\,\triangle\,}

\title{Слабо гомоклинические группы эргодических действий}
\author{В.В. Рыжиков}

\usepackage[T1]{fontenc} 
\usepackage[cp1251]{inputenc} 
\usepackage[russian]{babel}
\date{}
\begin{document}
\Large
\maketitle
\begin{abstract}{
 Гомоклиническая группа эргодического действия введенна М.И. Гординым. В статье рассмотены
естественные расширения этой группы, установлена    связь гомоклинических групп с факторами действия и  К-свойством. 
Показана траекторная полнота  слабо гомоклинических групп. Доказана
эргодичность слабо гомоклинических групп гауссовских и пуассоновских действий, установлена 
тривиальность гомоклинических групп для классов действий ранга 1.  Установлена  
связь слабо гомоклинических групп с такими  асимптотическими инвариантами  как жесткость действия, 
локльный ранг и слабое кратное перемешивание.  Обсуждаются  нерешенные задачи.     

Библиография: 25 названий. УДК: 517.987. MSC2010: 28D05 (Primary); 58F11 (Secondary).

Ключевые слова и фразы: эргодические действия, гомоклинические группы, гауссовские динамические системы, 
пуассоновские надстройки,   действия  ранга один, слабое кратное перемешивание.}
\end{abstract}

\section{Введение.  }
Пусть     $T$ --  автоморфизм стандартного вероятностного пространства $(X,\mu)$,  
 через $T$  также обозначим отвечающий ему унитарный оператор 
в  $L_2(X,\mu)$, заданный равенством $Tf(x)=f(Tx)$. 

{\bf    Гомоклиническая группа }  автоморфизма  $Т$  определена  М.И. Гординым  в \cite{G}  как
$$ H (T) = \{S\in Aut(\mu): T^{- n} ST^n \to I, \ n \to \infty \}, $$
где подразумевается  сильная операторная сходимость к тождественному оператору $I$.
Им была замечена связь свойств гомоклинической 
группы с перемешивающими свойствами автоморфизма  и инициировано  изучение  взаимосвязей 
со спектральными и энтропийными инвариантами  автоморфизма.
  
{\bf Взаимосвязи с инвариантами динамической системы.}  Эргодическая гомоклиническая группа легко обнаруживается у
бернуллиевских сдвигов. В связи с этим было интересно узнать, насколько тесно связано наличие такой группы с
 положительностью энтропии и лебеговским  спектром.     Дж.~Кинг, отвечая на вопрос Гордина, 
в  \cite{K} построил автоморфизм с нулевой энтропией 
и  эргодической гомоклинической  группой. Этот результат в некотором смысле был усилен автором: 
эргодической гомоклинической группой могут обладать автоморфизмы с простым сингулярным спектром. 
Нужные примеры были найдены в классе перемешивающих  
пуассоновских надстроек, все они обладают эргодической   гомоклинической группой  \cite{R}.  
Эргодичность гомоклинической группы является чисто метрическим инвариантом действия. Так, например,  
у геодезического потока эта группа эргодична, а у спектрально изоморфного ему орициклического потока она тривиальна 
(вытекает из результатов М.Ратнер \cite{R}).
Наличие эргодической гомоклинической группы автоморфизма влечет 
за собой свойство кратного перемешивания.  Это дало новое   доказательство кратного перемешивания 
 для перемешивающих пуассоновских надстроек, первоначально 
обнаруженное Э. Руа \cite{Roy}. Как мы увидим, причины, по которым перемешивающие   гауссовские динамические системы
обладают  кратным перемешиванием,  установленным В.П. Леоновым \cite{L},  также  имеют гомоклиническую
трактовку. 

Известное
$\Z^2$-действие  Ледрапье \cite{Le}  не обладает кратным перемешиванием, поэтому его 
нельзя  расширить до действия с эргодической гомоклинической группой. С.В. Тихонов  сообщил автору,
что из наблюдений в \cite{R} и  \cite{T}  вытекает тривиальность гомоклинической группы  этого действия.  
В настоящей статье доказано, что перемешивающие действия ранга 1 также имеют  тривальную 
слабо гомоклиническую группу.  

{\bf Слабо гомоклиническая и $P$-гомоклиническая группа.} Рассмотрим естественные расширения  
группы $H(T)$, содержательные для неперемешивающих систем и особенно для групповых действий. 
Слабо гомоклиническая  группа  $WH (T)$ определяется   как
$$ WH (T) = \{ S\in Aut(\mu): \frac {1} {N} \sum_{i = 1}^{N}T^{- i} ST^i \to I, 
\ N \to \infty \}. $$
Напомним, что здесь и ниже подразумевается сильная сходимость.  
Бесконечному множеству  $P\subset \Z$ и  $\Z$-действию, порожденному автоморфизмом $T$,
 отвечает группа
$$ H_P (T) = \{S\in Aut(\mu) \ : \ T^{- n} ST^n \to I, \ n \in P, \ n \to \infty \}. $$
Таким образом, автоморфизму сопоставляется континульное   семейство $P$-гомоклинических  групп. 
 Тривиальность, нетривиальность и 
эргодичность этих групп  являются метрическими   инвариантами  автоморфизма $T$. 
 В качестве примера, показывающего 
содержательность таких инвариантов, приведем утверждение, которое докажем позже.
\vspace {3mm}

\bf Теорема 1. \it Для любых бесконечных  множеств $C',D'\subset \N$  найдутся 
бесконечные  подмножества $C\subset C'$, $D\subset D'$  и слабо перемешивающие  
автоморфизмы $S,T$ такие, что
группы  $H_{C}(S)$, $H_{D}(T)$  эргодичны, а группы  $H_{C}(T)$, $H_{D}(S)$ тривиальны.
При этом $WH (S)$ и  $WH (T)$ эргодичны.
\rm
\vspace {3mm}

 Для  действия $\{T_g\}$ локально-компактной группы  $G$ и ее бесконечного подмножества  $P\subset G$ 
аналогичным образом определяется $P$-гомоклиническая группа. Заметим, что  $H(\{T_g\})=H_G(\{T_g\})$. 
Для действия Ледрапье  группа $H_{\Z^2}(\{T_g\})$  тривиальна, хотя $H_P(\{T_g\})$ эргодичны для  некоторых
бесконечных подмножеств $P$,  например, для одномерных подгрупп  в  $\Z^2$. 

 В статье установлены:

(\S 2)   траекторная полнота  $P$-гомоклинических групп, их связь с с факторами действия,   К-свойством, слабым кратным перемешиванием;
 
(\S 3)  эргодичность слабо гомоклинических групп гауссовских и пуассоновских действий;

(\S 4)  тривиальность слабо гомоклинических групп для классов действий ранга 1.

В \S  5  обсуждаются нерешенные задачи.   
 Неопределяемые понятия см., например, в \cite{KSF}.


\section{Полнота  гомоклинических групп,  факторы действия, \\  К-свойство и кратное перемешивание}

Полная группа $[S]$ автоморфизма $S$ по определению состоит из автоморфизмов   $R$ таких, что
для некоторой измеримой целочисленной функции $p(x)$ выполняется  $R(x)=S^{p(x)}$ 
(о траекторной теории и полных группах см. \cite{V85}). Полной группой $[\{S_g\}]$ 
группового действия $\{S_g\}$ называем группу, порожденную всеми полными 
группами $[S_g]$.
Если эргодическое действие полно, то оно содержит эргодический автоморфизм. Наиболее просто найти такой автоморфизм с 
рационально-дискретным спектром. По известной теореме Дая полная группа эргодического автоморфизма содержит 
 изоморфные копии  всех эргодических автоморфизмов. 
\vspace {3mm}

\bf Теорема 2.   \it  Группа   $H_P (\{T_g\})$ является полной группой: $$H_P (\{T_g\}) = 
[H_P (\{T_g\})].$$   Аналогичное утверждение верно для  группы $ WH(T)$.
\rm  
\vspace {3mm}

\bf Следствие. \it Если группы $H_P (T)$, $WH(T)$ эргодичны, то они 
с точностью до сопряжения содержат все эргодические автоморфизмы.
\rm
\vspace {3mm}
 
Доказательство теоремы  становится очевидным, если воспользоваться следующей интерпретацией:  
условие $T_{g(i)}^{-1} ST_{g(i)} \to I$  эквивалентно тому, что график функции $S:X\to X$ под действием
последовательности отображений $T_{g(i)}\times T_{g(i)}$ притягивается к диагонали в $X\times X$.  
Это следует понимать так:  для любого измеримого разбиения  $X=\bigsqcup_{k} E_k$  график функции $S$ под действием  
$T_{g(i)}\times T_{g(i)}$ асимптотически оказывается внутри $ D=\bigsqcup_{k} (E_k\times E_k)$.  Формально говоря, 
мера проекции той части графика функции $S$ на сомножитель $X$,  что  не попала в  $ D$,  стремится к 0.
Сказанное верно для всех степеней автоморфихма $S$.

Если график притягивается к диагонали, то части графика тоже  притягиваются.   Поэтому график, составленный 
из кусков графиков функций  $S^p$,  притягивающихся к диагонали, также будет притягиваться к диагонали. Но это 
эквивалентно утверждению теоремы  в случае групп $H_P (T)$.  Для  $WH(T)$ рассуждения аналогичны. Теорема
доказана.

Фактором преобразования $T$ называется его ограничение на инвариантную $\sigma$-алгебру.
\vspace {3mm}

\bf Теорема 3.  
 \it  Если  действие  $ \{T_g\} $ коммутативно и не имеет факторов,  то либо 
  группа $H_P (\{T_g\})$ эргодична, либо тривиальна: $ H_P (\{T_g\}) = \{I \} $.  Аналогичное утверждение верно 
для группы  $WH(T)$. \rm
\vspace {3mm}

Доказательство.  Пусть  для всех $S$  из  $H_P (\{T_g\})$ выполнено  $SA=A$.  Заметим, что $T_g^{-1}ST_g\in H_P (T)$.
Действительно, $$T_{g(i)}^{-1}T_g^{-1}ST_gT_{g(i)} = T_g^{-1} T_{g(i)}^{-1}ST_{g(i)}T_g\ \ \to \ I$$
(здесь мы воспользовались коммутативностью действия).  Тогда имеем  
$$T_g^{-1}ST_gA=A, \ \  ST_gA=T_gA.$$
  Значит, $\sigma$-алгебра неподвижных множеств относительно действия $H_P (T)$  является инвариантной 
относительно $T$.  Из этого непосредственно вытекает утверждение теоремы.
\vspace {3mm}

К-автоморфизмом является автоморфизм, не имеющий фактора с нулевой энтропией, 
по тереме Синая \cite{Si} он имеет бернуллиевский фактор. 
\vspace {3mm}

\bf Теорема 4.  
 \it Для К-автоморфизма $T$  группа  $H (T)$  является эргодической.\rm
\vspace {3mm}

Доказательство.  Пусть $H(T)=H_{\Z}(T)$  не является эргодической. Тогда алгебра $\F$ неподвижных относительно 
 $H (T)$ множеств
в силу доказанного  есть фактор автоморфизма  $T$.  Этот фактор нетривиальный, поэтому по тереме Синая он
имеет бернуллиевский подфактор $T'$, действующий на подадгебре $\F'\subset \F$. В силу теоремы Рохлина автоморфизм 
 $T$ можно представить как косое произведение над  $T'$,
 поэтому каждый  гомоклинический элемент $S'$ для  $T'$ продолжается до   
гомоклинического преобразования  $S\in H(T)$, где  косое произведение $S$ есть расширение 
 $S'$ с  тождественным преобразованием в слоях.  Но у бернуллиевского $T'$ группа $H (T')$ эргодична.
 Получили, что $\F'$ не содержит нетривиальных элементов, неподвижных относительно действия группы $H (T)$.
   Противоречие показывает, что алгебра $\F$ должна быть тривиальной, т.е.  $H (T)$ эргодична, что и требовалось.

\bf Контрпример Ф. Ледрапье.  \rm  Как мы знаем, эргодичность   гомоклинической группы действия влечет 
за собой свойство кратного перемешивания
\cite{R}. Рассмотрим перемешивающее $\Z^2$-действие Ледрапье
  $\{T_z\}$, для которого   выполняется 
$$\mu(A\cap T_{z_i}A\cap T_{w_i}A)=  0, \  {z_i},{w_i}, z_i-w_i \to\infty,$$
причем $$\mu(A)= \frac 1 2 , \ \  A  = T_{z_i}A\ \tr T_{w_i}A. $$
Покажем, следуя С.В. Тихонову, как отсюда вытекает тривиальность гомоклинической группы.

Пусть $S$  принадлежит гомоклинической группе действия $H(\{T_z\})$.  Тогда из леммы 4.2 \cite{R} получим
$$ \mu(SA\cap T_{z_i}A\cap T_{w_i}A)\to 0. $$
В  \cite{T} показано, что в этом случае 
$$ \mu(SA \tr T_{z_i}A  \tr T_{w_i}A)\to 0. $$
Следовательно, $\mu(SA \tr A)\to 0$, $SA=A$. Тогда  $ST_zA=T_zA$,  т.е.  $S$ является тождественным на 
всей алгебре, порожденной множествами  $T_zA$.  Таким образом, мы установили, что $S=I$, $H(\{T_z\})=\{I\}$. 

 Как устроена   слабо гомоклиническая группа  действия Ледрапье, мы не знаем, хотя
многие $H_P$ группы ($P\subset \Z^2$)  этого действия являются эргодическими.
\section{Слабо гомоклинические группы гауссовских  и пуассоновских действий}
В теории представлений   и эргодической теории хорошо известны 
гауссовские и пуассоновские действия (см. \cite{R},\cite{Roy},\cite{KSF}, \cite{V} --\cite{FR}). 
Первые возникают благодаря  представлению  ортогональной группы $O(\infty)$ в группе автоморфизмов пространства 
$\R^\infty$ с гауссовской мерой, 
а вторые (пуассоновские надстройки) -- действию группы автоморфизмов
 бесконечного пространства с мерой   на пространстве конфигураций с мерой Пуассона.
 Автоморфизму $T$  стандартного бесконечного пространства 
с мерой   соответствуют и пуассоновская    надстройка $P(T)$ и  гауссовскай автоморфизм  $G(T)$ 
(здесь $T$ рассматривается как
ортогональный оператор в соответствующем пространстве $L_2$).  
Автоморфизмы  $P(T)$ и    $G(T)$   унитарно эквивалентны  
оператору   $$exp(T)=\bigoplus_{n=0}^{\infty} T^{\odot n},$$ где  $T^{\odot 0}$ -- одномерный тождественный оператор, 
$T^{\odot n}$ -- симметрическая тензорная степень оператора $T$. При этом они могут обладать разными 
метрическими свойствами:   
 некоторые пуассоновские действия не имеют собственных факторов \cite{FR}, 
но все гауссовские системы обладают континуальной структурой факторов (см. \cite{RT}).  
\vspace {3mm}

\bf Теорема 5. \it Слабо перемешивающие гауссовские и пуассоновские системы обладают эргодической 
слабо гомоклинической группой.\rm
\vspace {3mm}

Доказательство. Сперва рассмотрим случай гауссовского действия. Обозначим через $F$ 
 группу ортогональных операторов в $ l_2$, которые сопряжены операторам вида $S=I + K$,
где $K$ -- конечномерный  оператор ($\| K\| \leq 2$). 

Пусть оператор $U\in O(\infty)$  имеет непрерывный спектр.
Покажем, что $F\subset WH(U)$.  Так как $U\otimes U$  также  имеет непрерывный спектр, при $\ N\to \infty$ 
выполняется 
$$ \frac {1} {N} \sum_{ 1}^{N} U^i\otimes U^i\ \to_w 0.$$
  Следовательно,   для любых векторов $f,g$ имеем 
$$ \frac {1} {N} \sum_{ 1}^{N} (U^i f,g)^2  \ \to 0.\eqno (1)$$
Пусть $\pi$ обозначает ортопроекцию на конечномерное пространство $Im K$ --  образ оператора $ K$. 
Выберем базис $\{g_k\}$  в $Im K$. В силу (1) выполняется 
$$ \frac {1} {N} \sum_{ 1}^{N} |(U^i f,g_k)|  \ \to 0,$$
что  с учетом  конечности  базиса  влечет за собой    
$$ \frac {1} {N} \sum_{ 1}^{N} \| \pi U^i f\|  \ \to 0.$$
Тогда
$$ \frac {1} {N} \sum_{ 1}^{N} \| K U^i f\| \leq  \frac {1} {N} \sum_{ 1}^{N} \|2 \pi U^i f\| \ \to 0, \ N\to \infty,$$
следовательно, 
$$\left\| f -  \frac {1} {N} \sum_{ 1}^{N}U^{- i} SU^i f \right\|= 
\left\| \frac {1} {N} \sum_{ 1}^{N}U^{- i} KU^i f \right\|  \leq
 \frac {1} {N} \sum_{ 1}^{N} \| KU^i f \| \ \to  0.$$
Этим мы установили сильную сходимость
$$  \frac {1} {N} \sum_{ 1}^{N}U^{- i} SU^i\ \to \ I.$$

Известно, что группа $F$ плотна в 
$O(\infty)$ в сильной операторной топологии (это вытекает, например,  из того факта, что любая матрица размера 
$n\times n$  вкладывается нужным образом в ортогональную матрицу размера $2n\times 2n$).  
Действие $G(O(\infty))$ эргодично, следовательно,  эргодично действие   $G(F)$.  Но  из $F\subset WH(U)$ вытекает, что 
$$G(F)\subset WH(G(U)).$$
Таким образом, слабо гомоклиническая группа гауссовского действия $G(U)$  содержит эргодическую подгруппу,
значит, она эргодична.

Теперь рассмотрим  пуассоновские действия. Теперь обозначим через $F$ группу всех финитных  автоморфизмов, т. е. автоморфизмов 
 $S$, совпадающих с 
тождественным автоморфизмом на бесконечном пространстве $(X,\mu)$ кроме множества $supp S$  конечной меры.
Пусть $T$  --  автоморфизм пространства  $(X,\mu)$ такой, что оператор $T$, ему соответствующий, 
 имеет непрерывный спектр. Сказанное  равносильно отсутствию у  $T$  собственного значения 1.
 В этом случае  имеем 
$$ \frac {1} {N} \sum_{ 1}^{N} (T^i f,g)  \ \to 0, \ N\to \infty.$$
Иначе ненулевой слабый предел для множества векторов $ \frac {1} {N} \sum_{ 1}^{N} T^if$ в силу усреднения
будет неподвижным вектором относительно $T$, что мы запрещаем.  Положив $f=\chi_A, g=\chi_B$  для множеств
$A,B$ конечной меры, получим
$$ \frac {1} {N} \sum_{ 1}^{N} \mu(T^i A\cap B)  \ \to 0, \ N\to \infty.$$
 Следовательно,
$$\left\| f -  \frac {1} {N} \sum_{ 1}^{N}T^{- i} ST^i f \right\|  \leq
 \frac {1} {N} \sum_{ 1}^{N} \| T^i f - ST^i f  \| \leq $$
$$\leq \frac {1} {N} \sum_{ 1}^{N} 2\mu (T^i A\cap supp S) 
\ \to  0, \ N\to \infty.$$
Получили
$$  \frac {1} {N} \sum_{ 1}^{N}U^{- i} SU^if\ \to \ f$$
для индикаторов $f$  множеств конечной меры, значит, это верно для всех $f\in L_2$.
Установлено, что 
$S\in  WH(T)$, а в силу непрерывности пуассоновского вложения $P$ имеем  $ P(S)\in WH(P(T)).$
Так как $F$ плотна в 
$Aut(X,\mu)$, а действие $P(Aut(X,\mu))$ эргодично, группа   $P(F)$ тоже  эргодична.
  Таким образом, слабо гомоклиническая группа $WH(P(T))$ пуассоновского действия $P(T)$  эргодична.
Теорема доказана.
\vspace {3mm}

Теперь мы можем доказать теорему 1.  Для этого над понадобится следующее утверждение.
\vspace {3mm}

\bf Теорема 6. \it Для любых бесконечных  множеств $C',D'\subset \N$  найдутся 
бесконечные  множества $C\subset C'$ и $D\subset D'$    и бесконечные эргодические автоморфизмы $\bar S,\bar T$ такие, что 

$\bar S^i  \to 0$ и $\bar T^i \to I$ при $i\to \infty,  i\in C$,
$\bar S^i  \to I$ и $\bar T^i \to 0$ при $i\to \infty,  i\in D$ \rm
\vspace {3mm}
\rm
\vspace {3mm}

Доказательство.  Cтроим конструкции $\bar S$   и $\bar T$  ранга 1  
(определение и обозначение 
параметров см. в \cite{R}) так,
чтобы высоты башен $h_j$ у них были одинаковыми (это условие не обязательно, но удобно). Полагаем $r_j\to\infty$.

На нечетном этапе $j$  у  конструкции $\bar T$   все надстройки, кроме последней,  
 имеют одинаковую высоту  $ s_j>h_j$, причем  $h_j+s_j\in C'$,
а у  конструкции $\bar S$ все надстройки, кроме последней,  
 пусть  имеют  высоту высоту $2s_j$.
Учитывая, что  $r_j\to\infty$, получим  $\bar T^{h_j+s_j}  \to I$  и   $\bar S^{h_j+s_j}  \to 0$,
если  $j$  пробегает нечетные числа. Множество значений $h_j+s_j$ для нечетных $j$  есть искомое  $C$.
 
Надстройки  над последними колоннами  выбираются так, чтобы высоты $h_{j+1}$ для  обеих конструкций совпадали. 
Например,  для конструкции $\bar S$ положим $s_j(r_j)=0$, а  для конструкции $\bar T$ пусть 
$s_j(r_j)=s_j(r_j-1)$.

На четном этапе $j$  у  конструкции $\bar S$   все надстройки, кроме последней,  
 имеют одинаковую высоту  $ s_j>h_j$,  $h_j+s_j\in D'$,
а у $\bar T$  все надстройки, кроме последней,  
 имеют  высоту высоту $ 2s_j$.
При $r_j\to\infty$ имеем   $\bar S^{h_j+s_j}  \to I$   $\bar T^{h_j+s_j}  \to 0$,
когда $j$  пробегает четные числа. Значения $h_j+s_j$ для четных $j$  образуют требуемое множество   $D$.
Надстройки  над последними колоннами  выбираются так, чтобы высоты $h_{j+1}$ у них совпадали.
Теорема доказана.

Доказательство теоремы 1.
Пользуясь теоремой 6, положим $S=P(\bar S)$, $T=P(\bar T)$.  Имеем $F\subset H_C(S)$ и $F\subset H_D(T)$,  следовательно,
группы  $H_{C}(S)$, $H_{D}(T)$  эргодичны. Группы  $H_{C}(T)$, $H_{D}(S)$ тривиальны, так как 
 $T\to I$ при $i\to \infty,  i\in C$ и  $S\to I$ при $i\to \infty,  i\in D$.  Таким образом, пуассоновские
надстройки (аналогично гауссовские) дают нужные примеры, теорема 1 доказана.

\bf Слабое кратное перемешивание. \rm
Чтобы упростить изложение,  будем рассматривать только перемешивание порядка 2.
Для более высоких порядков все совершенно аналогично.
Говорят, что слабо перемешивающий автоморфизм $G$ вероятностного пространства обладает 
слабым кратным перемешиванием  порядка 2, если  из  
$$  G^{m_i},G^{n_i},G^{n_i-m_i}\to\Theta  \eqno (2)$$
вытекает, что для любых измеримых множеств $A,B,C$  выполнено
$$ \mu(A\cap G^{m_i}B\cap G^{n_i}C)\to  \mu(A)\mu(B)\mu(C). \eqno (3)$$
Если действие перемешивающее и обладает слабым кратным перемешиванием, то оно обладает кратным перемешиванием.
\vspace {3mm}

\bf Теорема 7. \it Слабо перемешивающие гауссовские и пуассоновские системы обладают 
слабым кратным перемешиванием всех порядков.\rm
\vspace {3mm}

Доказательство.  Пусть для гауссовского автоморфизма $G=G(U)$ 
выполнено (2). Тогда имеем
$  U^{m_i},U^{n_i}\to 0.$
Для любого автоморфизма $S$ имеем
$$ \mu(A\cap G^{m_i}B\cap G^{n_i}C)=\mu(SA\cap SG^{m_i}B\cap SG^{n_i}C)=$$
$$
=\mu(SA\cap G^{m_i}G^{-m_i}SG^{m_i}B\cap G^{n_i}G^{-n_i}SG^{n_i}C).$$
Пусть теперь  $S\in G(F)$, где  $F$ -- группа операторов, рассмотренная в доказательстве теоремы 4. 
Так как $G(F)\subset H_M(T)\cap H_N(T)$,  получим
$$
G^{-m_i}SG^{m_i}, \ G^{-n_i}SG^{n_i} \to \ I, $$
$$ \mu(A\cap G^{m_i}B\cap G^{n_i}C)-\mu(SA\cap G^{m_i}B\cap G^{n_i}C)\to 0. $$
Из этого  с учетом того, что 
действие $G(F)$ эргодично, а из  $G^{n_i-m_i}\to\Theta$
вытекает 
$$\mu(G^{m_i}B\cap G^{n_i}C)\to \mu(B)\mu(C), $$ 
следует (3),  см.  \cite{R}, лемма 4.2.

Для пуассоновских надстроек рассуждаем аналогично:  вместо оператора  $U$ рассматриваем автоморфизм 
$T$ бесконечного пространства с мерой.
Отвечающая ему пуассоновская надстройка $P=P(T)$ обладает теми же  свойствами, что и гауссовский автоморфизм $G$.
 Роль группы $F$ теперь
выполняет группа всех финитных  автоморфизмов бесконечного пространства,  см.  доказательство теоремы 5.
Таким образом, слабо перемешивающие  пуассоновские действия  обладают слабым кратным перемешиванием по тем же причинам, 
что и гауссовские: все  гомоклинические группы $H_{\{m_i\}}$ при условии, что   последовательность $\{m_i\}$ 
перемешивающая, содержат общую эргодическую подгруппу (надcтройку над группой $F$). Но это обеспечивает 
свойство кратного перемешивания.  Теорема доказана.

\section{Ранг действия и   слабо гомоклиническая группа}

Сохраняющее меру обратимое 
преобразование
   $T: X \to X$ пространства Лебега $ (X, \mu) $ имеет ранг 1, если 
 существует последовательность $$ \xi_j= \{E_j, \ TE_j, \ T^2 E_j, \ \dots, T^{ h_j-1} E_j, \tilde {E}_j \}  $$
 измеримых разбиений  пространства 
 $ X $ таких, что 
 любое множество конечной меры  
 с заданной точностью аппроксимируется 
$\xi_j$-измеримыми множествами  для всех больших 
  значений $j$.
Это понятие  родственно циклической аппроксимации (без скорости аппроксимации), см. \cite{KS}.
Перемешивающие  конструкции ранга 1 появились в \cite{O}, они оказались полезны 
для построения разнообразных контрпримеров в эргодической теории \cite{Ru}.
Техника  аппроксимаций самоприсоединений (self-joinings) почти инвариантными мерами  позволила 
 получить новые результаты о действиях конечного ранга, см. \cite{R93}.
  Как мы увидим,  почти инвариантные меры  весьма удобны для исследования гомоклинических групп  действий ранга 1. 
\vspace {3mm}

\bf Теорема 8. \it Если перемешивающий автоморфизм $T$  вероятностного пространства имеет ранг 1, 
то его слабо гомоклиническая группа  тривиальна:  $WH(T)=\{I\}$.  \rm
\vspace {3mm}

Доказательство.  Пусть $S\in WH(T)$.    Положим  
$$U^a_j= \bigcup_{i=1}^{[ah_j]} T^{i} E_j\times T^{i}. $$
Наша цель --  установить, что $S=Id$. Для этого будет показано, что часть графика автоморфизма $S$, 
которая находится  внутри множества $(U^a_j\times U^a_j)$, 
 $\frac 1 2 < a<1$,  
асимптотически сосредоточена в околодиагональном множестве
$$\bigcup_{i=1}^{[ah_j]} (T^iE_j\times T^iE_j).$$
Но это возможно только в том случае, когда график автоморфизма $S$ является  диагональю в $X\times X$, 
что нам и нужно.

Фиксируем $a$,    $1<2a<2$.   Переходя при необходимости к подпоследовательности по $j$ (обозначим ее снова через $j$) 
 имеем
 $$\mu(S U^a_j\cap U^a_j)\to c\geq 2a-1>0.$$
Обозначим через $\Delta_S$ меру, для  которой выполнено $\Delta_S(A\times B)=\mu(SA\cap B)$  для всех измеримых $A,B$.
Пусть $$\sigma_j=\Delta_S(\ \ \ | U^a_j\times U^a_j),$$ положим $M_j=[(1-a)h_j]$
и
$$ \eta_j=\frac 1 {M_j} \sum_{i=1}^{M_j} (T^{i}\times T^{i}) \sigma_j$$
Тогда $ \eta_j\to  \Delta=\Delta_{Id},$
так как предельная точка для последовательности мер  $\eta_j$  является  $T\times T$-инвариантной мерой, 
абсолютно непрерывной относительно предела
мер
$$\frac 1 {M_j}\sum_{i=1}^{M_j} (T^{i}\times T^{i}) \Delta_S, $$
который равен  $\Delta$,   так как $S\in WH(T)$.

Представим 
$$ \sigma_j =\sum_{p,q=1}^{M_j}\sigma_j^{p,q},$$
где  меры $ \sigma_j^{p,q}$ определяются равенством  
$$ \sigma_j^{p,q}(Y) = \sigma_j(Y\cap (T^{p}E_j\times T^{q}E_j)).$$
Тогда
$$\eta_j= \sum_{p,q=1}^{M_j} \eta_j^{p,q},$$
где
$$\eta_j^{p,q}=\frac 1 {M_j} \sum_{i=1}^{M_j}(T^{i}\times T^{i}) \sigma_j^{p,q}.$$

Если $A,B$ --   $ \xi_j$-измеримые множества, то 
$$\eta_j^{p,q}(A\times B)=c^{p,q}_j\Delta^{p,q}_j(A\times B),$$
где $$c^{p,q}_j=\mu(ST^pE_j\cap T^qE_j)/\mu(SU^a_j\cap U^a_j),$$
а $\Delta^{p,q}_j$  -- нормированная часть меры $\Delta_{T^{p-q}}$,
расположенная в объединении $\cup_{i=1}^{M_j} (T^{i+p}E_j\times T^{i+q}E_j)$.
Итак, 
$$\eta_j(A\times B)=\sum_{p,q=1}^{M_j}c^{p,q}_j\Delta^{p,q}_j(A\times B).$$
$$\sum_{p,q=1}^{M_j}c^{p,q}_j\Delta^{p,q}_j \ \to \Delta.$$
Но почти инвариантные меры  $\Delta^{p,q}_j$ близки к $\Delta$  только при $p=q$, так 
как в силу свойства перемешивания при   большом $M$ и 
$|p-q| \geq M$ меры $\Delta_{T^{p-q}}$ и их почти инвариантные части близки к мере $\mu\times\mu$.
А при $0\neq |m| < M$ меры $\Delta_{T^m}$ и с ними $\Delta^{p,p+m}_j$  отделены от меры $\Delta$.
Таким образом, для любого положительного числа $a<1$   меры 
$(U^a_j\times U^a_j)\Delta_S$  
асимптотически сосредоточены в околодиагональных множествах
$$\bigcup_{i=1}^{[ah_j]} (T^iE_j\times T^iE_j).$$
Это возможно, только в случае, когда  $\Delta_S=\Delta$, т.е.  $S=Id$.  Теорема доказана.

В конце доказательства теоремы 8 мы использовали свойство перемешивания.  Что меняется, если его не предполагать?
Мы получим, что части  $\Delta^{p(j),q(j)}_j$ сдвигов диагональной меры близки к $\Delta$, поэтому  
 имеет место следующее утверждение:

\it Если  автоморфизм $T$  вероятностного пространства имеет ранг 1 и его  
слабо гомоклиническая группа  нетривиальна, то он обладает свойством частичной жесткости с показателем $\frac 1 2$:
т. е.  для некоторой последовательности $m(i)\to\infty$ и некоторого марковского оператора $P$ выполняется
 $$T^{m(i)}\to \frac 1 2 I +\frac 1 2 P.$$ \rm
\vspace {3mm}

Причина, по которой получается коэффициент $\frac 1 2$, состоит в том, что условие 
$\Delta_S(U^a_j\times U^a_j)>c>0$ гарантированно выполняется при $a> \frac 1 2$.
Но, если предположить, что $S$ -- апериодический автоморфизм, для любого $a>0$  найдется $k>0$  такое, что
для бесконечного набора $j$  выполнено
$$\Delta_{S^k}(U^a_j\times U^a_j)>c>0.$$
Из этого мы получаем, что $T$   обладает свойством частичной жесткости с показателем $1-a$, но в силу произвольности
$a$ получаем показатель 1, т.е. $T$ является жестким автоморфизмом.
Тем самым доказано следующее утверждение.
\vspace {3mm}
 
\bf Теорема 9.  \it Если  автоморфизм  вероятностного пространства $T$  имеет ранг 1 и его  
слабо гомоклиническая группа $WH(T)$ содержит апериодический автоморфизм, то $T$ является жестким:
$T^{m(j)}\to I$  для некоторой последовательности $m(j)\to\infty$. \rm
\vspace {3mm}

Таким образом,  ограниченные нежесткие конструкции ранга 1  (см. \cite{R13}), не обладают  
эргодической слабо гомоклинической группой.
Как следствие мы получили    гомоклиническое доказательство  утверждения из 
 \cite{FR}:  \it если пуассоновская надстройка обладает рангом 1, то она жесткая.  \rm

\bf  Локальный ранг. \rm Укажем на некоторую  связь слабо гомоклинических групп с локальным рангом 
(об этом инварианте см., например, \cite{R1}).
 Пусть $T$ - автоморфизм стандартного пространства Лебега $ (X, \mu) $, $ \mu (X) = 1 $.
Пусть для числа $ \beta> 0 $
существует последовательность $ \xi_j $ разбиений $ X $ вида
$$ \xi_j = \{B_j, TB_j, T^2B_j, \dots, T^{h_j-1} B_j, C_j^1, \dots, C_j^{m_j} \dots \},$$
стремящаяся к разбиению на точки, назовем ее аппроксимирующей, причем 
 $ \mu (U_j) \to \beta $, где
$ U_j = \bigsqcup_ {0 \leq k <h_j} T^kB_j $.

Локальный ранг $ \beta (T) $ определяется как максимальное число $ \beta $
для которого автоморфизм $ T $ обладает соответствующей последовательностью аппроксимирующих разбиений.
\vspace {3mm}

 \bf Теорема 10.\it 
Если эргодический автоморфизм $T$ не является частично жестким с показателем жесткости $\beta(T)$ и    $ \beta (T) > 0 $,
то группа $ WH (T) $ тривиальна.
\rm

Доказательство этого утверждения  мало отличается от доказательства теоремы 9, оставляем его в качестве упражнения.

\section{ Заключительные замечания и вопросы.}

Сформулируем несколько   вопросов, возникающих в связи  с полученными результатами.

1. Понятие  $P$-гомоклинических групп без  изменений переносится в топологическую динамику, где
 хорошо известны  аналоги свойств эргодичности,  ранга,  перемешивания, кратного перемешивания, энтропии. 
Контрпримеры к топологическому кратному перемешиванию (см. \cite{DK}, \cite{GM}) показывают,  что
гомоклинические группы соответствующих  действий не обладают свойством минимальности (аналогом эргодичности).
 Свойства  гомоклинических групп  могут представить некоторый  интерес в изучении
топологических инвариантов непрерывных отображений.   Верно ли, что гомеоморфизм компакта с положительной 
топологической энтропией обладает нетривиальной гомоклинической группой?

 2. Какие  структуры гомоклинических групп  возможны? Например, бывают ли нетривиальные гомоклинические группы, 
являющиеся полной  группой  периодического преобразования? 
  Есть ли слабо перемешивающий автоморфизм, у которого все $P$-группы тривиальны? Мы знаем, что таким свойством
обладают  автоморфизмы  с дискретным спектром. 

 3.  Если гомоклиническая  группа  одного действия  тривиальна, а другого зргодична, 
будут ли они дизъюнктны:    
 марковский оператор, сплетающий эти действия, совпадает с ортопроекцией на пространство констант?
Как связана гомоклиническая группа действия с гомоклиническими группами ее факторов?
Бернулиевские факторы могут порождать действие с неэргодической гомоклинической группой. 

 4. Верно ли, что  действия с бедной структурой самоприсоединений (например, с  минимальными самоприсоединениями,
 см. \cite{Ru})
имеют тривиальную слабо  гомоклиническую группу?  В случае перемешивающих действий ранга 1
это так.  
Для неперемешивающих систем  также имеются соответствующие примеры:
вполне эргодические нежесткие ограниченные  конструкции  ранга 1 обладают минимальными самоприсоединениями 
(см. \cite{R13}), а  их  слабо гомоклиническая группа тривиальна в силу теоремы 7.
Потоками бесконечного ранга с тривиальной слабо гомоклинической  группой являются унипотентные потоки 
(следствие результатов  \cite{Ra}).

 5. Что можно сказать про гомоклинические группы  слабо перемешивающего жесткого автоморфизма $T$  ранга 1  или 
положительного  локального ранга? Может ли некоторая группа $H_P(T)$ быть эргодической?
Если  слабо гомоклиническая группа такого автоморфизма всегда неэргодична, то  это дало бы  другое  
доказательство теоремы \cite{Rue} о том, что гауссовские автоморфизмы не обладают рангом 1, 
и аналогичный, но   новый результат  про пуассоновские  надстройки.

Московский государственный университет

E-mail: vryz@mail.ru

\end{document}